\documentclass[11pt]{amsart}

\usepackage[margin=1.1in]{geometry}
\usepackage[utf8]{inputenc}

\usepackage[romanian]{babel}
\usepackage{combelow}

\usepackage{amsmath}
\usepackage{amsthm}
\usepackage{amssymb}
\usepackage{colonequals}
\usepackage{tikz}
\usepackage{hyperref}
\usepackage{xypic}

\newcommand{\PP}{\mathbb{P}}

\renewcommand{\tilde}{\widetilde}
\renewcommand{\hat}{\widehat}
\renewcommand{\phi}{\varphi}
\renewcommand{\setminus}{\smallsetminus}

\newtheorem{thm}{Theorem}[section]

\newtheorem{cor}[thm]{Corollary}
\newtheorem{prop}[thm]{Proposition}

\theoremstyle{definition}

\theoremstyle{remark}
\newtheorem{rem}[thm]{Remark}

\def\OO{\mathcal{O}}

\def\cM{\mathcal{M}}

\def\cU{\mathcal{U}}

\def\mm{\overline{\mathcal{M}}}

\newcommand\dJ{de Jonqui\`{e}res }

\def\cc{\overline{\mathcal{C}}}

\title{Generalized de Jonqui\`{e}res divisors on generic curves}

\author[G. Farkas]{Gavril Farkas}
\address{Humboldt-Universit\"at zu Berlin, Institut F\"ur Mathematik,  Unter den Linden 6
\hfill \newline\texttt{}
 \indent 10099 Berlin, Germany} \email{{\tt farkas@math.hu-berlin.de}}

\dedicatory{To the memory of Csaba Varga (1959-2021)}

\begin{document}

\begin{abstract}
The classical \dJ  and MacDonald formulas describe the virtual number of divisors with prescribed multiplicities in a linear system on an algebraic curve. We discuss the enumerative validity of the \dJ formulas for a general curve of genus $g$.
\end{abstract}

\maketitle

\section{Introduction}

De Jonqui\`{e}res' formula \cite{dJ} is concerned with the following classical enumerative question: Given a suitably general (singular) plane curve of $\Gamma$ degree $d$ and geometric genus $g$, how many plane curves of given degree meet  $\Gamma$ in $n_i$ unspecified points with contact order $a_i$, for $i=1, \ldots, e$? De Jonqui\`{e}res using an ingenious recursive argument (later considerably simplified by Torelli \cite{T} and then slightly generalized by Allen \cite{All}) showed that the number in question equals
$$
\frac{\bigl[a_1^{n_1} a_2^{n_2}\cdots a_e^{n_e}\bigr]}{n_1!\ n_2! \ \cdots n_e!}, \ \ \  \  \mbox{ where in general we define the quantity}$$

\begin{equation}\label{eq:dJ1}
[a_1 \cdots a_e]=a_1\cdots a_e\ \frac{g!}{(g-e-1)!}\left(\frac{a_1\cdots a_e}{g-e}-\sum_{i=1}^e \frac{a_1 \cdots \hat{a_i}\cdots a_e}{g-e+1}+\cdots +(-1)^e \frac{1}{g}\right).
\end{equation}

The formula (\ref{eq:dJ1}) recovers many well known formulas in the theory of algebraic curves, for instance the number $2^{g-1}(2^g-1)$ of odd theta characteristics on a smooth curve of genus $g$, or the Pl\"ucker formula for the total number of ramification points on a linear series on a curve. The original proofs \cite{dJ}, \cite{T} of the \dJ formula use an induction on the multiplicities $a_i$ coupled with the \emph{Brill-Cayley correspondence principle}. For a historic perspective on the \dJ formula we refer to Zeuthen's treatise \cite[136]{Z}, or if one prefers English, the books of Coolidge  \cite[Book 3, Chapter 3.3]{Coo} or Baker \cite[pages 35-45]{Ba}.  De Jonqui\`{e}res' formula has been rediscovered by MacDonald \cite{McD} and Vainsencher \cite{V} and a summary of their work, reinterpreting this number as a fundamental class of a modified diagonal on the symmetric product of a smooth curve can be found in the book \cite{ACGH}.

\vskip 4pt

In order to formulate the problem in modern terms, let $C$ be a smooth curve of genus $g$ and we fix a linear series $\ell=(L,V)\in G^r_d(C)$. For a partition $\mu=(a_1, \ldots, a_e)$ of $d$, we define the \dJ cycle $DJ_{\mu}(C, \ell)$ to be the locus of divisors of the type $a_1\cdot x_1+\cdots+a_e\cdot x_e$
lying in the linear system $\ell$. Observe that $DJ_{\mu}(C, \ell)$ can be realized as the rank $r$ degeneracy locus of the evaluation morphism of vector bundles
$$\chi\colon V\otimes \OO_{C^e}\longrightarrow J_{\mu}(L)$$ over the product $C^e$, where the fibre of the vector bundle $J_{\mu}(L)$ over a point $(x_1, \ldots, x_e)$ equals the $d$-dimensional vector space $L_{| a_1\cdot x_1+\cdots+a_e\cdot x_e}$.  Accordingly, the virtual dimension of $DJ_{\mu}(C, \ell)$ equals $e-d+r$. In the case $e=d-r$, this number equals zero and one expects $\ell$ to contain finitely many divisors with multiplicities prescribed by the partition $\mu$. As pointed out in \cite[page 359]{ACGH}, the virtual class of this degeneracy locus can be realized via the Porteous formula as the coefficient of the monomial $t_1 t_2 \cdots t_e$ in the polynomial
$$\bigl(1+a_1^2 t_1+\cdots+a_e^2 t_e\bigr)^g \bigl(1+a_1 t_1+\cdots +a_e t_e\bigr)^{d-r-g}.$$
It is straightforward to see that this is simply a convenient way to repackage compactly the information contained in the formula (\ref{eq:dJ1}).
For instance, we obtain that a linear system $\ell \in G^r_d(C)$ is expected to contain precisely
$$2^r\left( {d-r\choose r}+ g{d-r-1\choose r-1}+ {g\choose 2}{d-r-2\choose r-2}+\cdots\right)$$
divisors containing $r$ double points, that is, of the type $2\cdot x_1+\cdots+2\cdot x_r+x_{r+1}+\cdots+x_{d-r}$ and so on. Here we use the convention that ${m\choose -h}=0$ when $h>0$.

\vskip 5pt

More generally, we consider a positive partition $\mu=(a_1, \ldots, a_e)$ and  set $|\mu|:=a_1+\cdots+a_e$ and $\ell(\mu):=e$. For $0\leq f\leq |\mu|$ we define the generalized \dJ (secant) locus
$$DJ_{\mu}^f(C, \ell):=Z_{|\mu|-f}(\chi)=\Bigl\{(x_1, \ldots, x_e)\in C^e:\mathrm{dim}\ \bigl|V(-a_1\cdot x_1-\cdots-a_e\cdot x_e)\bigr|\geq r-|\mu|+f\Bigr\}.$$
Being a degeneracy locus, each  component of $DJ_{\mu}^f(C, \ell)$ has dimension at least $e-f(r+1-|\mu|+f)$. If $\mu=(1^e)$, then using the notations of \cite{CM} or \cite{F}, we observe that $DJ_{\mu}^f(\ell)=V_e^{e-f}(\ell)$  can be identified with the variety of $e$-secant $(e-f-1)$-planes to the embedded curve $C\stackrel{|V|}\hookrightarrow \PP^r$. Moreover, if $|\mu|=d$, then $DJ_{\mu}^{d-r}(C, \ell)=DJ_{\mu}(C, \ell)$ is the locus of \dJ divisors in the linear series $\ell$. De Jonqui\`{e}res loci have been used to study the geometry of the moduli spaces of curves or that of strata of holomorphic differentials \cite{BCGGM}. The class of effective divisors on $\mm_g$ involving \dJ conditions have been computed in \cite{Cot},  \cite{FV}, \cite{FV2}, or \cite{Mu}.

\vskip 4pt

The question of how to interpret the \dJ count when a curve $C\subseteq \PP^r$ acquires singularities has been treated both in classical and modern times. The problem we address in this note on the other hand is the enumerative validity of the \dJ count when $C$ is a general curve in moduli. We treat this problem variationally and consider  \dJ cycles  associated to all linear systems $\ell\in G^r_d(C)$, that is, we set up the correspondence:
\begin{equation}\label{eq:corr}
\xymatrix{
  & \Sigma_{\mu}^f(C):=\Bigl\{\bigl(\ell, x_1, \ldots, x_e\bigr): (x_1, \ldots, x_e)\in DJ_{\mu}^f(C, \ell)\Bigr\} \ar[dl]_{\pi_1} \ar[dr]^{\pi_2} & \\
   G^r_d(C) & & C^e      \\
                 }
\end{equation}

The main result of this paper is then summarized as follows:
\begin{thm}\label{thm:main}
Let $C$ be a general curve of genus $g$ and we fix a partition $\mu=(a_1, \ldots, a_e)$, as well as positive integers $d, r$ and $f$  with $\rho(g, r,d)\geq 0$ and $|\mu|-r\leq f\leq |\mu|$. Then
each irreducible component of $\Sigma_{\mu}^f(C)$ has dimension $\rho(g,r,d)+e-f(r+1-|\mu|+f)$.
Accordingly, if $$\rho(g,r,d)+e-f(r+1-|\mu|+f)<0,$$ then $DJ_{\mu}^f(C, \ell)=\emptyset$ for every linear series $\ell\in G^r_d(C)$.
\end{thm}

\vskip 3pt

This result generalizes \cite[Theorem 0.1]{F} to the case of an arbitrary partition $\mu$, the result in \emph{loc.cit.} corresponding to the case when $\mu=(1^e)$.  It also generalizes Ungureanu's results \cite[Theorem 1.5]{U} corresponding to the case when $|\mu|=d=\mbox{deg}(\ell)$, asserting that if $C$ is a general curve, no linear series $\ell\in G^r_d(C)$ possesses a \dJ divisor of length $e<d-r$. Observe that the case $f=|\mu|-r$ in Theorem \ref{thm:main} can be obviously reduced to the classical \dJ case, by extending the partition $\mu$  to $\mu'=(\mu, 1^{d-|\mu|})$ of the degree $d$ of the curve in question.

\vskip 4pt

We now discuss several cases in which Theorem \ref{thm:main} applies. The first case beyond the classical \dJ situation treated for instance (under some restrictive assumptions) in \cite{U} is when $f=|\mu|+1-r$, when the residual linear series $\bigl|V(-a_1\cdot x_1-\cdots -a_e\cdot x_e)\bigr|$ is a pencil, which can be formulated as saying that under the map
$\phi_{\ell}\colon C\rightarrow \PP^r$ induced by the linear series $\ell$, the $(a_i-1)$-st osculating planes to $C$ at the points $x_i$ span a codimension two plane, that is,
\begin{equation}\label{eq:osc}
\bigl\langle a_1\cdot x_1, \ldots, a_e\cdot x_e\bigr\rangle \cong \PP^{r-2}.
\end{equation}

\vskip 5pt

\noindent {\bf{Tangential secants.}}
Let us consider the case $a_1=2$ and $a_2=\cdots=a_e=1$ and $f=1$, in which case  the condition (\ref{eq:osc}) translates into saying that
$\langle 2\cdot x_1, x_2, \ldots, x_e\rangle \cong \PP^{e-1}$, that is, the tangent line to $C$ at the point $x_1$ lies in the $(e-1)$-plane spanned by the points $x_1, \ldots, x_e$. Following classical terminology, we say that $\langle x_1, \ldots, x_e\rangle$ is a \emph{tangential $(e+1)$-secant} to $C$. Theorem \ref{thm:main} can be formulated in this case as follows:

\begin{cor}\label{cor:tang}
We fix positive integer $g, r, d$ and $e$ such that $2e< r+1-\rho(g,r,d)$. For a general curve $C$ of genus $g$, no linear seris $\ell\in G^r_d(C)$ carries a tangential $(e+1)$-secant.
\end{cor}

Note that every space curve $C\subseteq \PP^3$ of degree $d$ and genus $g$ is expected to have finitely many tangential \emph{trisecants} and their number $T(d,g)=2(d-2)(d-3)+2g(d-6)$, which can derived from the \dJ formula,  has been first computed by Salmon and Zeuthen \cite[64]{Z}, see also \cite[page 364]{ACGH}. It is an interesting result of Kaji \cite{K}, valid to a large extent even in positive characteristic, that an \emph{arbitrary} smooth space curve $C\subseteq \PP^3$ cannot have infinitely many tangential trisecants, see also \cite{BP} for various extensions of this result. For space curves, our Corollary \ref{cor:tang} reduces to the Brill-Noether Theorem, but already for curves $C\subseteq \PP^4$ it goes beyond that and it states  that when $\rho(g, r,d)=0$ a general such curve has no tangential trisecants.

\vskip 5pt

\noindent {\bf{Multiple tangents.}} Passing now to the case of tangent planes, that is, when $a_1=\cdots =a_e=2$, we look at $(2e-2)$-planes in $\PP^r$ that are tangent to $C$ at $e$ points, that is,
$$\bigl\langle 2\cdot x_1, \ldots, 2\cdot x_e\bigr\rangle\cong \PP^{2e-2}.$$
We call such a configuration an \emph{degenerate $e$-tangent} to $C\subseteq \PP^r$. With this terminology, Theorem \ref{thm:main} takes the following form:

\begin{cor}\label{cor:tg2}
Fix positive integers $g,r,d, e$ with $\rho(g, r, d)\geq 0$ and $3e<r+2-\rho(g,r,d)$. Then a general curve $C$ of genus $g$ has no linear series $\ell\in G^r_d(C)$ with degenerate $e$-tangents.
\end{cor}

The simplest case where Corollary \ref{cor:tg2} applies is when $e=2, r=5$. It says that for a general curve $C$ of genus $g$, no embedded curve
$\varphi_{\ell}\colon C\rightarrow \PP^5$ of degree $d$ with $\rho(g, r, d)=0$ has a pair of coplanar tangent lines.

\vskip 3pt

Another immediate application of Theorem \ref{thm:main} is when again $a_1=\cdots=a_e=2$ but this time $f=2e-r>0$, hence
$$\bigl \langle 2\cdot x_1, \ldots, 2\cdot x_e\bigl \rangle\cong \PP^{r-1}.$$
In other words, the points $x_1, \ldots, x_e$ span a \emph{tangent hyperplane}. We find the following result:

\begin{cor}\label{cor:tg3}
Fix integers $g\geq 1$, $r\geq 3$ and $d$ such that $\rho(g, r,d)\geq 0$ and $e\geq r+1$. Then for a general curve $C$ of genus $g$ the locus of linear systems $\ell\in G^r_d(C)$ such that $\varphi_{\ell}\colon C\hookrightarrow \PP^r$ admits an $e$-secant tangent hyperplane is equal to $\rho(g,r,d)+r-e$.
\end{cor}

In particular, for $e=r+1$ specializes to the known result \cite{U}, that for a Brill-Noether general curve $C\subseteq \PP^r$ no hyperplane can be tangent at more than $r$ points.

\vskip 5pt

\noindent {\bf{Flex lines and bitangents.}}
A general smooth space curve $C\subseteq \PP^3$ is expected to possess no bitangent or flex lines lines, that is, no \dJ divisors of length two corresponding to the partitions $\mu=(2, 2)$ and $\mu=(3, 1)$ respectively. We consider the problem more generally for curves $C\subseteq \PP^r$ and our result in this case lends a sharp  form to this expectation.

\begin{cor}\label{cor:flex}
Fix positive integers $g\geq 1$, $r\geq 3$ and $d$ with $\rho(g, r,d)\geq 0$ and $a_1, a_2$ such that $$a_1+a_2>\frac{\rho(g,r,d)+2r}{r-1}.$$
Then for a general curve $C$ of genus $g$, no degree $d$ embedding
$\varphi_{\ell}\colon C\hookrightarrow \PP^r$ possesses a secant line meeting the image of $C$ with multiplicities $a_1$ and $a_2$ at the points of secancy.
\end{cor}

\vskip 3pt

For instance when $r=3$, $e=2$ and $|\mu|=4$, Corollary \ref{cor:flex} implies that when $\rho(g, 3,d)\leq 1$, for a general curve $C$ of genus $g$ no embedding $\varphi_{\ell}\colon C\hookrightarrow \PP^3$ of degree $d$ possesses either a bitangent or a flex line.

\vskip 5pt

The last application of Theorem \ref{thm:main} is to the case when the partition $\mu$ is of length one.

\begin{cor}\label{cor:1}
We fix positive integers $g, r, d$ and $a$ such that $2a>\rho(g,r,d)-1+2r$. Then a general curve $C$ of genus $g$ carries no linear series $\ell\in G^r_d(C)$ having a point $x\in C$ with $\ell(-a\cdot x)\in G^1_{d-a}(C)$.
\end{cor}

Specializing even further to the case $d=2g-2$ and $r=g-1$ in which case $\ell$ necessarily equals the canonical linear series $|\omega_C|$, via the Riemann-Roch Theorem Corollary \ref{cor:1} can be reformulated as stating that for a general curve of genus $g$, if $a\geq g-1$ we have that
$$h^0\bigl(C, \OO_C(a\cdot x)\bigr)\leq a+2-g,$$
for each point $x\in C$. When $a=g-1$ we obtain that $C$ carries no pencil of degree $g-1$ totally ramified at a point, which is a well-known result.  The locus of curves $[C]\in \cM_g$ having such a pencil has been studied by Diaz \cite{D}, who also computed the class of its compactification in $\mm_g$.

\vskip 8pt

\noindent  {\small{{\bf{Acknowledgments}}:
This article is dedicated to the memory of Csaba Varga (1959-2021), one of the defining mathematical personalities in Cluj/Kolozsv\'ar. The intellectual influence he had on my mathematical development during my studies at the Babe\cb{s}-Bolyai University between 1991 and 1995  cannot be overstated.

\vskip 4pt

The author was supported by the DFG Grant \emph{Syzygien und Moduli} and by the ERC Advanced Grant SYZYGY of  the European Research Council (ERC) under the European Union Horizon 2020 research and innovation program (grant agreement No. 834172).
}}

\section{Generalized \dJ divisors on flag curves}

We fix a smooth curve $C$ of genus $g$ and we denote by $G^r_d(C)$ the variety of linear systems of type $g^r_d$ on $C$, that is, pairs $\ell=(L, V)$, where $L\in \mbox{Pic}^d(C)$ and $V\subseteq H^0(C,L)$ is an $(r+1)$-dimensional subspace of sections. Recall that when $C$ is a general curve of genus $g$, then $G^r_d(C)$ is a smooth variety of dimension equal to the \emph{Brill-Noether number} $\rho(g,r,d)=g-(r+1)(g-d+r)$. Our proof of Theorem \ref{thm:main} is by degeneration and we will use throughout the theory of limit linear series. We begin by quickly recalling the notation for vanishing and ramification sequences of linear series on curves largely following \cite{EH1} and \cite{EH2}.

\vskip 4pt

If $\ell=(L,V)\in G^r_d(C)$ is a linear series, the \emph{ramification sequence} of $\ell$ at a point $q\in C$
\[
\alpha^{\ell}(q) : 0\leq \alpha_0^{\ell}(q) \leq \cdots \leq \alpha_r^{\ell}(q)\leq d-r
\]
is obtained from the \emph{vanishing sequence}
\[
a^{\ell}(q) : 0\leq a_0^{\ell}(q) < \cdots < a_r^{\ell}(q) \leq d
\]
by setting $\alpha^{\ell}_i(q) := a^{\ell}_i(q) - i$, for $i=0, \ldots, r$.  In case the underlying line bundle $L$ is clear from the context, we write
$\alpha^{V}(q) = \alpha^{\ell}(q)$ and $a^{V}(q) = a^{\ell}(q)$.  The \emph{ramification weight} of $q$ with respect to $\ell$ is defined as the quantity $\mathrm{wt}^{\ell}(q) := \sum_{i=0}^r \alpha^{\ell}_i(q)$.  We denote by
$$\rho(\ell, q) := \rho(g,r,d) - \mathrm{wt}^{\ell}(q)$$ the \emph{adjusted Brill-Noether number} of $\ell$ with respect to $q$.
We recall also the \emph{Pl\"ucker formula}
\begin{equation}\label{eq:pluecker}
\sum_{q\in C} \alpha^{\ell}(q)=(r+1)d+(r+1)r(g-1),
\end{equation}
measuring the total ramification of $\ell$. Incidentally, assuming that $\ell$ has only simple ramification points, that is, points with ramification sequence at most $(0,\ldots, 0,1)$, then (\ref{eq:pluecker}) is an instance of the \dJ formula (\ref{eq:dJ1}) applied to the linear series $\ell$ and to  the partition $\mu=(r+1, 1^{d-r-1})$ of $d$).

\vskip 4pt

Following Eisenbud-Harris \cite[page 364]{EH1}, let us recall that a \emph{limit linear series} on a curve $X$  of compact type consists of a collection $\ell=\bigl\{(L_C, V_C) \in G^r_d(C): C\mbox{ is a component of } X\bigr\}$, satisfying a  compatibility condition on the vanishing sequences at the nodes of $X$ in terms of the vanishing sequences of the aspects on the two (smooth) components of $X$ on which each node of $X$ lies. We denote by $\overline{G}^r_d(X)$ the variety of limit linear series of type $g^r_d$ on $X$. More generally, if $q\in X_{\mathrm{req}}$ is a smooth point and $\alpha=\bigl(0\leq \alpha_0\leq \cdots \leq \alpha_r\leq d-r\bigr)$ is a \emph{Schubert index}, we denote by $\overline{G}^r_d\bigl(X, (q,\alpha)\bigr)$ the variety of limit linear series $\ell\in \overline{G}^r_d(X)$ satisfying the condition $\alpha^{\ell}(q)\geq \alpha$. From basic principles it follows that each component has dimension at least $\rho\bigl(g, r, d, \alpha)=\rho(g, r,d)-\mathrm{wt}(q)$. Eisenbud and Harris offer in \cite[Theorem 1.1]{EH2} sufficient conditions ensuring when the equality
\begin{equation}\label{eq:strongBN}
\mbox{dim } \overline{G}^r_d\bigl(X,(q, \alpha)\bigr)=\rho(g,r,d)-\mathrm{wt}(\alpha)
\end{equation} holds, which we will make an essential use of in the course of proving Theorem \ref{thm:main}. In case a pointed curve $[X, q]$ satisfies the condition (\ref{eq:strongBN}) for each $r, d\geq 1$ such that $\rho(g, r,d)\geq 0$ and for each choice of a Schubert index $\alpha$, we say that $[X,q]$ verifies the strong Brill-Noether Theorem.

\vskip 4pt

Having fixed a positive partition $\mu=(a_1, \ldots, a_e)$, a positive integer $f$ with  $|\mu|-r\leq f\leq |\mu|$ and a smooth curve $C$, we have defined in the Introduction the subvariety $\Sigma_{\mu}^f(C)\subseteq G^r_d(C)\times C^e$. Due to its determinantal structure, each irreducible component of $\Sigma_{\mu}^f(C)$ has dimension at least
$$\mbox{dim } G^r_d(C)+e-f(r+1-|\mu|+f)\geq \rho(g,r,d)+e-f(r+1-|\mu|+f).$$
From this fact we obtain  that once one shows that for a general curve $C$ of genus $g$ each irreducible component of $\Sigma_{\mu}^f(C)$ has dimension \emph{at most} $\rho(g,r,d)+e-f(r+1-|\mu|+f)$, it will also follow that $\Sigma_{\mu}^f(C)$ is in fact equidimensional of this dimension.

\vskip 4pt

Assume we are in a situation when $\Sigma_{\mu}^f(C)$ in nonempty for a general (and therefore for an arbitrary) smooth curve $C$ of genus $g$.

\vskip 4pt

\subsection{Universal \dJ divisors on curves of compact type.} The proof of Theorem \ref{thm:main} relies, like several other proofs involving limit linear series, on degenerating a smooth curve of genus $g$ to a flag curve consisting of a rational spine and $g$ smooth elliptic tails. It is known \cite{EH1} and \cite{EH2} that such curves satisfy the Brill-Noether Theorem \emph{independently} of the position of the $g$ points of attachment on the rational spine. One has however to deal with the serious complication that, under this degeneration, although one has a good understanding of the aspects of the limit linear series on the flag curve, a priori there is no control on the position of the $e$ marked points lying in the support of a generalized \dJ divisor. For the combinatorial argument required to prove Theorem \ref{thm:main} it is however essential to ensure that one can always find such a flag curve degeneration of a generic curve of genus $g$ in which these $e$ marked points specialize to a subcurve of the flag curve having relatively small arithmetic genus. To make sure this is possible, we employ a strategy already used in \cite{F}, which relies on considering \emph{all} flag curves of genus $g$ at once and using certain basic facts about the geometry of the (rational) parameter space of such a curves.

\vskip 4pt

We  set some further notation. Let $j\colon \mm_{0,
g}\rightarrow \mm_g$ the map assigning to a
stable rational pointed curve $[R, p_1, \ldots, p_g]\in \mm_{0, g}$ \emph{fixed} smooth elliptic
tails $E_1, \ldots, E_g$ at the marked points $p_1, \ldots, p_g$. We denote the resulting compact type curve by
$$X:=R\cup_{p_1} E_1\cup \ldots \cup_{p_g} E_g,$$ that is, $p_a(X)=g$ and
let $p_R\colon X\rightarrow R$ be the map contracting each elliptic component $E_i$ to the point $p_i$.
We introduce the universal $n$-pointed curve $\cc_{g, n}=\mm_{g,n+1}$ of genus $g$ and denote by
$\pi\colon \cc_{g, n}\rightarrow \mm_{g, n}$ the morphism forgetting the
$(n+1)$-st marked point. For $e\geq 1$, we write $\pi_e\colon \cc_{g,n}^e \rightarrow
\mm_{g, n}$ for the $e$-fold fibre product of $\cc_{g, n}$ over
$\mm_{g, n}$.  We finally introduce the map
\begin{equation}\label{eq:mapchi}
\chi\colon \mm_{0, g}\times_{\mm_g} \cc_g^e\rightarrow \cc_{0, g}^e,
\end{equation}
which collapses the fixed elliptic \emph{tails} $E_1, \ldots, E_g$ and projects the corresponding marked points onto the rational \emph{spine} $R$. With the notation introduced above, we thus have
$$\chi\Bigl([R, p_1,\ldots, p_g], (x_1, \ldots, x_e)\Bigr)=\Bigl([R, p_1,
\ldots, p_g], \ p_R(x_1), \ldots,  p_R(x_e)\Bigr),$$
where $x_1, \ldots, x_e\in X$.

\vskip 4pt

Let $\overline{\mathfrak{DJ}} \subseteq \cc_g^e$ be the closure of the
locus of generalized \dJ divisors on smooth curves of genus $g$, that is, of the following determinantal variety

\begin{align*}
\mathfrak{DJ}:=\Bigl\{[C, x_1, \ldots, x_e]: [C]\in \cM_g, \ x_i\in C,  \ \exists \ell=(L, V) \in G^r_d(C) \mbox{
such that } \\ \mbox{ dim } 
\bigl|V(-a_1\cdot x_1-\cdots -a_e\cdot x_e)\bigr|\geq r-|\mu|+f\Bigr\}.
\end{align*}
Since we assume that $\Sigma_{\mu}^f(C)\neq \emptyset$ for a general curve $[C]\in \cM_g$, we have that $\pi_e\bigl(\overline{\mathfrak{DJ}}\bigr)=\mm_g$, where recall that $\pi_e\colon \cc_g^e\rightarrow \mm_g$.  Next, we define the locus
\begin{equation}\label{eq:defU}
\cU:=\chi\Bigl(\mm_{0,g}\times_{\mm_g}   \overline{\mathfrak{DJ}}\Bigr)\subseteq \cc_{0,g}^e.
\end{equation}
We use the commutativity of the following diagram, where the horizontal upper arrow is induced via the \emph{stabilization} isomorphism $\cc_{g,n}\cong \mm_{g,n+1}$, see \cite[page 175]{Kn} by taking fibre products
\[
\xymatrix{
 \cc_{0,g}^e \ar[d]^{\pi_e} \ar[r] & \cc_g^e \ar[d]^{\pi_e} \\
 \mm_{0,g} \ar[r]^{j} &  \mm_g \\
}
\]
in order to conclude that $\pi_e(\cU)=\mm_{0,g}$. We denote by $e-m$ the generic fibre dimension of the map
$\pi_{e | \cU}\colon  \cU \rightarrow \mm_{0, g}$. Thus $0\leq m\leq e$ and
$$\mbox{dim}\bigl(\cU\cap \pi_e^{-1}[R, p_1, \ldots, p_g]\bigr)=e-m,$$ for a general stable curve
$[R, p_1, \ldots, p_g]\in \mm_{0, g}$.

We introduce the birational map
$$\vartheta \colon \cc_{0, g}^e\rightarrow \mm_{0,
4}^{g-3+e}\cong (\PP^1)^{g-3+e}$$ whose components are the forgetful morphisms $\pi_i\colon \mm_{0, g+e}\rightarrow \mm_{0,4}$ which for $i=4, \ldots, g+e$ only retain the marked points labelled by $1, 2, 3$ and $i$ respectively.  Fixing for instance the first three marked points as usual $p_1=0$, $p_2=1$ and $p_3=\infty \in \PP^1$, by slightly abusing notation we can think of $\vartheta$ as the map assigning
$$\bigl([R, p_1, \ldots, p_g], x_1, \ldots, x_e\bigr)\stackrel{\vartheta}\mapsto \bigl(p_4, \ldots, p_g, x_1, \ldots, x_e\bigr)\in (\PP^1)^{g-3+e}.$$

Using essentially only the elementary fact that the diagonal of $\PP^1\times \PP^1$ is ample, we then establish in \cite[Proposition 2.2]{F}, that depending on whether  $\vartheta(\cU)\subseteq (\PP^1)^{g-3+e}$ intersects the small diagonal $(x_1=\cdots=x_e)$ in $(\PP^1)^{g-3+e}$ or not, one of the following three possibilities occur:

\vskip 5pt

\noindent $\bullet$ There exists a point $(p_4, \ldots, p_g, x_1, \ldots, x_e)\in \vartheta(\cU)$ with $x_1=\cdots=x_e$ and at least $g-m-3$ of the points $p_4, \ldots, p_g$ are mutually distinct.

\noindent $\bullet$ There exists a point $(p_4, \ldots, p_g, x_1, \ldots, x_e)\in \vartheta(\cU)$ such that at least $g-m$ of the points $p_4, \ldots, p_g$ are equal to a point $r\in \PP^1\setminus \{x_1, \ldots, x_e\}$.

\noindent $\bullet$ There exists a point $(p_4, \ldots, p_g, x_1, \ldots, x_e)\in \vartheta(\cU)$ such that $e-1$ of the marked points $x_1, \ldots, x_e$ are equal and at least $g-m$ of the points $p_4, \ldots, p_g$ are equal to $0$.

\vskip 4pt

Investigating the fibres of the map $\vartheta$ in each of these cases we find the following, see \cite{F}:

\begin{prop}\label{schub}
Keeping the notation above, if $\mathrm{dim}(\cU)=g-3+e-m$, there exists a point
$$\bigl([R, p_1, \ldots, p_g], x_1, \ldots,
x_e\bigr) \in \mm_{0,g}\times _{\mm_g} \overline{\mathfrak{DJ}},$$ such that on the flag curve
$X=R\cup_{p_1} E_1\cup \ldots \cup_{p_g} E_g$ the limiting \dJ divisor $(x_1, \ldots, x_e)$ satisfies either  (i)  $x_1=\cdots =x_e\in R\setminus \{p_1, \ldots,
p_g\}$, or else, (ii) $x_1, \ldots, x_e$ all lie on a connected subcurve $Y\subseteq X$ of genus at  most $m$ and with $\bigl|Y\cap
(\overline{X\setminus Y})\bigr|\leq 1$.
\end{prop}

\vskip 5pt

\subsection{The proof of Theorem \ref{thm:main}}\label{subsect:proof}
We fix a partition $\mu=(a_1, \ldots, a_e)$  and a positive integer $f\geq |\mu|-r$. We assume
that the variety $\Sigma_{\mu}^f(C)\subseteq G^r_d(C)\times C^e$ is not empty for every smooth curve $C$ of genus $g$. Keeping the notation above, we denote by $e-m$ the fibre dimension of the surjective morphism $\pi_e\colon \cU\rightarrow \mm_{0,e}$. Recall that we defined $\overline{\mathfrak{DJ}}\subseteq \cc_g^e$ to be the closure of the universal locus of \dJ divisors and we assume that $e-n$ is the generic fibre dimension of the surjective morphism
$$\pi_{e|\overline{\mathfrak{DJ}}}\colon \overline{\mathfrak{DJ}} \rightarrow \mm_g.$$
Since when specializing to the subvariety of flag curves via the map $j\colon \mm_{0,g}\hookrightarrow \mm_g$ the fibre dimension of $\pi_e$ can only go up, we have that $m\leq n$.   We now apply Proposition \ref{schub} and let
$X=R\cup E_1\cup \ldots \cup E_g$
be the corresponding flag curve of genus $g$ as above, where for $i=1, \ldots, g$ we denote by $p_i\in R$ the node corresponding to the intersection of the \emph{spine} $R$ (which may itself well be reducible) with the subtree of $X$ ending in the elliptic \emph{tail} $E_i$.  We
denote by $Y\subseteq X$ the connected subcurve of $X$ onto which the
marked points $x_1, \ldots, x_e$ (limiting a generalized \dJ divisor) specialize. According to Proposition \ref{schub} there are two  possibilities:

\vskip 5pt

\noindent {\bf{(i)}} $p_a(Y)=m\leq \mathrm{min}\{e, g\}$, or

\noindent {\bf{(ii)}} $x_1=\cdots =x_e\in R\setminus \{p_1, \ldots, p_g\}$.

\vskip 5pt

We first treat case (i). Let $Y':=\overline{X\setminus Y}$ be the subcurve of $X$ complementary to $Y$ and set
$\{p\}:=Y\cap Y'$.  When $m=g$, then set $Y:= X$ and $Y'=\emptyset$ and we let $p\in
X$ be a general (smooth) point.  The divisor $a_1\cdot x_1+\cdots+a_e\cdot x_e$ is a limit of generalized \dJ divisors on smooth curves of genus $g$ neighboring the genus $g$ curve of compact type $X$. Applying the formalism of stable reduction, we can find a flat family of nodal curves of genus $g$
$$\phi\colon \mathcal{X}\rightarrow (T, t_0)$$ over a smooth pointed curve, together with sections $s_1, \ldots, s_e\colon T\rightarrow \mathcal{X}$ such that:

\vskip 5pt

\noindent ${{\bf{(1)}}}$  The generic fibre $\phi^{-1}(t)=X_t$ is a smooth curve of genus $g$, whereas the central fibre  $$\tilde{X}:=\phi^{-1}(0)$$ is stably equivalent to $X$, that is, it is a curve of arithmetic genus $g$ obtained from $X$ by possibly attaching chains of smooth rational curves at the singularities of $X$.

\noindent ${{\bf(2)}}$ $s_i(0)=x_i\in \tilde{X}_{\mathrm{reg}}$ for
all $i=1, \ldots, e$.

\noindent ${{\bf{(3)}}}$  There exists a line bundle $L_{\eta}$ of related degree $d$ defined on the complement of the central fibre $X_{\eta}=\mathcal{X}\setminus \phi^{-1}(0)$,  and a subvector bundle  $V_{\eta}\subseteq \phi_*(L_{\eta})$ of rank $r+1$, such that for $t\neq 0$, setting $L_t=L_{\eta|X_t}\in \mbox{Pic}(X_t)$ and $V_t=V_{\eta|t}\subseteq H^0(X_t, L_t)$, we have
that $$\Bigl((L_t, V_t), s_1(0), \ldots, s_e(t)\Bigr)\in \Sigma_{\mu}^f(X_t),$$
that is, $\mbox{dim} \ \bigl|V_t(-a_1\cdot s_1(t)-\cdots-a_e\cdot s_e(t))\bigr|\geq r-|\mu|+f$.

\vskip 4pt

We shall denote by $\tilde{Y}\subseteq \tilde{X}$ the inverse image of $Y$ under the contraction morphism $\tilde{X}\rightarrow X$. Then set $\tilde{Y}':=\overline{\tilde{X}\setminus \tilde{Y}}$ and we still denote by $p$ the point of intersection of $\tilde{Y}$ and $\tilde{Y}'$.

Since when forming the family $\mathcal{X}\rightarrow T$ we allow us the possibility of a further base change and that of resolving the resulting singularities,  we may furthermore assume that the flag curve $\tilde{X}$ carries a (refined) limit linear series
$$\ell=\Bigl\{\ell_Z=(L_Z,V_Z): Z \mbox{ is a component of } \tilde{X}\Bigr\}\in \overline{G}^r_d\bigl(\tilde{X}\bigr)$$
obtained following the procedure described by Eisenbud and Harris \cite{EH1} as a limit of the linear series $(L_t, V_t)$.
Furthermore, the sublinear series described in (3) induce a limit linear series
$$\ell'=\Bigl\{\ell_Z'=(L_Z(-D_Z), V_Z'): Z \mbox{ is a component of }  \tilde{X}\Bigr\}\in \overline{G}^{r-|\mu|+f}_{d-|\mu|}\bigl(\tilde{X}\bigr),$$
where $D_Z$ is an effective divisor on $Z$ supported on the union of the points $s_1(0), \ldots, s_e(0)$ that happen to lie on $Z$ and the point
of intersection $Z\cap \overline{\tilde{X}\setminus Z}$ (which is a smooth point of $Z$), and $V_Z'\subseteq H^0(Z, L_Z')$ is respectively a subspace of sections of dimension
$r+1-|\mu|+f$.

\vskip 5pt

Note that $p$ is a smooth point of both subcurves $\tilde{Y}$ and $\tilde{Y}'$ of $\tilde{X}$, therefore it is a smooth point of a unique irreducible component of $\tilde{Y}$, respectively of a unique irreducible component of $\tilde{Y}'$. We consider the respective aspects of $\ell$ and slightly abusing notation, we denote by
$$a^{\ell_{\tilde{Y}}}(p)=\bigl(a_0< \cdots <a_r\bigr)$$ the sequence obtained by ordering the vanishing orders at $p$ of the sections corresponding to the irreducible component of $\tilde{Y}$ containing $p$.  Similarly,  we let
$$a^{\ell_{\tilde{Y}'}}(p)=\bigl(b_0< \cdots <b_r\bigr)$$ be the sequence obtained by ordering the vanishing orders at $p$ of the sections contained in the aspect of $\ell$ corresponding to the component of $\tilde{Y}'$ containing $p$.  Note that $a_i+b_{r-i}=d$ for $i=0, \ldots, r$.
Furthermore, by ordering the vanishing orders at $p$ of the aspect of $\ell'$ corresponding to the component of $\tilde{Y}$ containing $p$,  we obtain the sequence
$$a^{\ell'_{\tilde{Y}}}(p)=\bigl(a_{i_0}<\cdots
<a_{i_{r-|\mu|+f}}\bigr).$$
Clearly, this is a subsequence of $a^{\ell_{\tilde{Y}}}(p)$. The entries in the complementary subsequence can be ordered as well and we denote this subsequence by $\bigl(a_{j_0}<a_{j_1}< \cdots <a_{j_{|\mu|-f-1}}\bigr)$. Note that
$$\bigl\{a_{i_0}, \ldots, a_{i_{r-|\mu|+f}}\bigr\}\cup \bigl\{a_{j_0}, \ldots, a_{j_{|\mu|-f-1}}\bigr\}=\bigl\{a_0, \ldots, a_r\bigr\}.$$

\vskip 4pt

While the entries in the sequence $(a_{j_0}< \cdots <a_{j_{|\mu|-f-1}})$ corresponding to vanishing orders of sections of a linear series on a single irreducible component  of $\tilde{Y}$, using the procedure described in  \cite[Lemma 2.1]{F}, one can construct a sublimit linear series $\ell_{\tilde{Y}}^{\sharp}\in \overline{G}_d^{|\mu|-f-1}(\tilde{Y})$ of $\ell_{\tilde{Y}}$  such that its vanishing sequence $a^{\ell_{\tilde{Y}}^{\sharp}}(p)$ equals precisely $\bigl(a_{j_0}<\cdots < a_{|\mu|-f-1}\bigr).$

\vskip 4pt

We first assume  $\tilde{Y}'\neq
\emptyset$. The point $p\in \tilde{Y}$ is a smooth point and lies on one of its rational component. In particular the genus $m$ pointed curve $[\tilde{Y}, p]$ verifies the strong Brill-Noether Theorem, that is, both varieties $\overline{G}^{r-|\mu|+f}_{d-|\mu|}\bigl(\tilde{Y}, (p, \alpha^{\ell_{\tilde{Y}}'}(p))\bigr)$ and $\overline{G}_d^{|\mu|-f-1}\bigl(\tilde{Y}, (p, \alpha^{\ell_{\tilde{Y}}^{\sharp}(p)})\bigr)$ have the expected dimension given by the corresponding  adjusted Brill-Noether numbers, in particular these numbers must be non-negative, cf.
\cite[Theorem 1.1]{EH2}. We thus obtain the following two inequalities by writing this for out for the limit linear series $\ell_{\tilde{Y}}$ and $\ell_{\tilde{Y}'}$ respectively:
\begin{equation}\label{egy1}
\mbox{dim } \overline{G}_d^{|\mu|-f-1}\Bigl(\tilde{Y}, \bigl(p, \alpha^{\ell_{\tilde{Y}}^{\sharp}(p)}\bigr)\Bigr)=\rho\bigl(\ell_{\tilde{Y}}^{\sharp}, p\bigr)=\rho\bigl(m, |\mu|-f-1, d\bigr)-a_{j_0}-\cdots-a_{j_{|\mu|-f-1}}+{|\mu|-f\choose 2}\geq
0,
\end{equation}
as well as
\begin{equation}\label{egy2}
\mbox{dim } \overline{G}^{r-|\mu|+f}_{d-|\mu|}\Bigl(\tilde{Y}, \bigl(p, \alpha^{\ell_{\tilde{Y}}'}(p)\bigr)\Bigr)=\rho(\ell'_{\tilde{Y}}, p)=
\end{equation}
$$\rho\bigl(m, r-|\mu|+f,
d-|\mu|\bigr)-a_{i_0}-\cdots- a_{r-|\mu|+f}+{r+1-|\mu|+f\choose 2}\geq 0.
$$

The same considerations can be applied to the complementary subcurve $\tilde{Y}'$ of $\tilde{X}$. The point of attachment $p$ lies on a rational component component of $\tilde{Y}'$, therefore the strong Brill-Noether inequality holds for $\ell_{Y'}$ as well, and we obtain:

\begin{equation}\label{egy3}
\mbox{dim }\overline{G}^r_d\Bigl(\tilde{Y}', (p, \alpha^{\ell_{\tilde{Y}'}}(p))\Bigr)=\rho(\ell_{\tilde{Y}'}, p)=\rho(g-m, r, d)-\bigl(b_0+\cdots+b_r\bigr)+{r+1\choose 2}\geq 0.
 \end{equation}

 We add the inequalities (\ref{egy1}), (\ref{egy2}) and (\ref{egy3}) together and
 use the fact that $(\ell_{\tilde{Y}}, \ell_{\tilde{Y}'})$ form a refined limit linear series, therefore
 the vanishing orders of $\ell_{\tilde{Y}}'$, $\ell_{\tilde{Y}}^{\sharp}$ and those of $\ell_{\tilde{Y}'}$ respectively add up, that is,
 $$\sum_{k=0}^r b_k+\sum_{k=0}^{r-|\mu|+f}
 a_{i_k}+\sum_{k=0}^{|\mu|-f-1} a_{j_k}=\sum_{k=0}^r \bigl(a_k+b_{r-k}\bigr)=(r+1)d.$$

 We obtain the following estimate:

 \begin{align*}\label{eq:ali1}
 0\leq \rho(g-m, r,d)+\rho\bigl(m, r-|\mu|+f, d-|\mu|\bigr)+\rho\bigl(m, |\mu|-f-1,d\bigr)\\
 -(r+1)d+{r+1\choose 2}+{r+1-|\mu|+f\choose 2}+{|\mu|-f\choose 2}\\
 =\rho(g,r,d)-f(r+1-|\mu|+f)+m \leq \rho(g, r,d)-f(r+1-|\mu|+f)+e,
 \end{align*}
which is precisely the second half of Theorem \ref{thm:main}. Note that in the last inequality, the assumption
$m\leq e$ guaranteed by Proposition \ref{schub} is absolutely essential.

\vskip 5pt

In the case $m=g$, when necessarily $e\geq g$ and $\tilde{Y}=\tilde{X}$, we proceed along similar lines. We add together
inequalities (\ref{egy1}) and (\ref{egy2}) to obtain:

 \begin{align*}
 \rho(g, r, d)+e-f(r+1-|\mu|+f)=
 \left(\rho\bigl(g, r-|\mu|+f,
 d-|\mu|\bigr)-\sum_{k=0}^{r-|\mu|+f} a_{i_k}+{r+1-|\mu|+f\choose
 2}\right)\\
 +\left(\rho\bigl(g, |\mu|-f-1, d\bigr)-\sum_{k=0}^{e-f-1}a_{j_k}+{e-f\choose
 2}\right)+\sum_{k=0}^{r-|\mu|+f} a_{i_k}+\sum_{k=0}^{|\mu|-f-1}
 a_{j_k}-{r+1\choose 2}+e-g\\
 =\mbox{dim } \overline{G}^{r-|\mu|+f}_{d-|\mu|}\Bigl(\tilde{Y}, \bigl(p, \alpha^{\ell_{\tilde{Y}}'}(p)\bigr)\Bigr)+\mbox{dim } \overline{G}_d^{|\mu|-f-1}\Bigl(\tilde{Y}, \bigl(p, \alpha^{\ell_{\tilde{Y}}^{\sharp}(p)}\bigr)\Bigr)\\+\sum_{k=0}^{r-|\mu|+f} a_{i_k}+\sum_{k=0}^{|\mu|-f-1}
 a_{j_k}-{r+1\choose 2}+e-g\geq 0,
 \end{align*}
 since $$\sum_{k=0}^{r-|\mu|+f} a_{i_k}+\sum_{k=0}^{|\mu|-f-1} a_{j_k}=\sum_{k=0}^r a_k\geq {r+1\choose 2}$$ and, as explained, $e\geq g$.

\vskip 5pt

 Assume finally we are in the case (ii), that is, when $x_1=\cdots =x_e\in
 R\setminus \{p_1, \ldots, p_g\}$. Keeping the previous notation, we observe that the limit linear series $\ell \in \overline{G}^r_d(\tilde{X})$ has vanishing sequence
 at $x_1$
$$a^{\ell}(x_1)\geq \Bigl(0, 1, \ldots, |\mu|-f-1, |\mu|, |\mu|+1, \ldots, r+f-1, r+f\Bigr),$$ therefore 
$\mathrm{wt}^{\ell}(x_1)\geq f(r+1-|\mu|+f)$. Taking into account that $[\tilde{X}, q]$ satisfies the strong Brill-Noether Theorem, cf.
 \cite[Theorem 1.1]{EH2}, Theorem 1.1, we obtain the inequality
$$0\leq \mbox{dim }\overline{G}^r_d\Bigl(\tilde{X}, \bigl(x_1, \alpha^{\ell}(x_1)\Bigr)\leq \rho(g, r, d)-f(r+1-|\mu|+f)\leq \rho(g, r, d)+e-f(r+1-|\mu|+f).$$

This concludes the proof that the assumption $\Sigma_{\mu}^f(C)\neq \emptyset$ for a general curve of genus $g$ implies that $\rho(g,r,d)+e-f(r+1-|\mu|+f)\geq 0$. 
 
\vskip 5pt

We come now to the dimensionality statement for the variety $\Sigma_{\mu}^f(C)\subseteq G^r_d(C)\times C^e$, when $C$ is a general curve of genus $g$. Recalling from the Introduction that $\pi_2\colon \Sigma_{\mu}^f(C)\rightarrow C^e$ is the natural projection, with our notation we have $\mbox{dim } \pi_2\bigl(\Sigma^f_{\mu}(C)\bigr)=e-n\leq e-m$, where $e-n$ has been defined as the minimal fibre dimension of the surjection $\overline{\mathfrak{DJ}}\rightarrow \mm_g$. We now estimate the fibre dimension of $\pi_2$ over a general point $(y_1, \ldots, y_e)\in \pi_2(\Sigma^f_{\mu})$.  To that end, we specialize once more to the locus of flag curves. For an $e$-pointed curve $[X, x_1, \ldots, x_e]$ of compact type, where the marked points are pairwise distinct smooth points of $X$, we denote by $\Sigma_{\mu}^f(X, x_1, \ldots, x_e)$ the subvariety of $\overline{G}^r_d(X)$ consisting of limit linear series
$$\ell=\Bigl\{\ell_Z=(\ell_Z,V_Z): Z \mbox{ is a component of }X\Bigr\}\in \overline{G}^r_d(X)$$
possessing a sublimit linear series of the form
$$\ell'=\Bigl\{\ell_Z'=\bigl(L_Z(-D_Z), V_Z'\bigr): Z \mbox{ is a component of }X\Bigr\}\in \overline{G}^{r-|\mu|+f}_{d-|\mu|}(X),$$
where $\mbox{supp}(D_Z)=Z\cap \bigl(\overline{(X\setminus Z)}\cup\{x_1, \ldots, x_e\}\bigr)$. As already explained, via Proposition \ref{schub} we may consider a further degeneration to a flag curve  $[\tilde{X}, x_1, \ldots, x_e]$, where $\tilde{X}=\tilde{Y}\cup \tilde{Y}'$ with $\tilde{Y}\cap \tilde{Y}'=\{p\}$ satisfies the conditions (1)-(3). Recall that $x_1, \ldots, x_e \in \tilde{Y}_{\mathrm{req}}\setminus \{p\}$.  It follows that for the generic fibre dimension of $\pi_2$ the following inequality holds:

$$\mbox{dim }\pi_2^{-1}(y_1, \ldots, y_e)\leq \mbox{dim }\Sigma_{\mu}^f\bigl(\tilde{X}, x_1, \ldots, x_e\bigr).$$
Furthermore, the dimension of $\Sigma_{\mu}^f\bigl(\tilde{X}, x_1, \ldots, x_e)$ cannot exceed the dimension of the space of triples $\bigl(\ell'_{\tilde{Y}}, \ell^{\sharp}_{\tilde{Y}}, \ell_{\tilde{Y}'}\bigr)$ described earlier, which as explained, via the estimates (\ref{egy1}), (\ref{egy2}) and  (\ref{egy3}) equals

\begin{align*}
\mbox{dim } \overline{G}^{r-|\mu|+f}_{d-|\mu|}\Bigl(\tilde{Y}, \bigl(p, \alpha^{\ell_{\tilde{Y}}'}(p)\bigr)\Bigr)+\mbox{dim } \overline{G}_d^{|\mu|-f-1}\Bigl(\tilde{Y}, \bigl(p, \alpha^{\ell_{\tilde{Y}}^{\sharp}(p)}\bigr)\Bigr)+\mbox{dim }\overline{G}^r_d\Bigl(\tilde{Y}', (p, \alpha^{\ell_{\tilde{Y}'}}(p))\Bigr)\\
=\rho(g,r,d)-f(r+1-|\mu|+f)+m.
\end{align*}
It follows that
\begin{align*}
\mbox{dim } \Sigma^f_{\mu}(C)\leq \mbox{dim } \pi_2\bigl(\Sigma_{\mu}^f(C)\bigr)+\mbox{dim }\Sigma_{\mu}^f\bigl(\tilde{X}, x_1, \ldots, x_e\bigr)\\
\leq e-n+m+\rho(g,r,d)-f(r+1-|\mu|+f)\leq e-f(r+1-|\mu|+f),
\end{align*}
since, as explained,  $m\leq n$. This  brings the proof of Theorem \ref{thm:main} to an end.

$\hfill$ $\Box$

\begin{rem}
A natural extension of Theorem \ref{thm:main} could be to consider the transversality of curves $C\subseteq \PP^r$ with respect to non-linear spaces. For instance, staying at the level of space curves, it is expected that a general curve $C\subseteq \PP^3$ has finitely many $8$-secant conics (but no $9$-secant conics), finitely many $12$-secant twisted cubics (but not $13$-secant twisted cubics) and so on. The smooth curves confirming this expectation have been recently characterized as those for which the blow-up of $\PP^r$ along $C$ yields a threefold with big and nef anticanonical divisor, see \cite{BL}. The (virtual) number of $8$-secant conics to $C\subseteq \PP^3$ has been computed by Katz \cite{K} as an iteration of multiple point formulas. It would be interesting to have a study of the enumerative validity of this and other similar formulas mirroring Theorem \ref{thm:main}. In this case however more subtle phenomena, related to the (Strong) Maximal Rank Conjecture \cite[Conjecture 5.1]{AF}, must come into play and which go beyond the Brill-Noether genericity of the curve in question. It is for instance clear that whenever $C\subseteq \PP^3$ lies on a quadric there is a positive dimensional family of $8$-secant conics, so at the very least these curves will have to be excluded, probably other as well.
\end{rem}

\bibliographystyle{plain}

\end{document}